\documentclass[conference]{IEEEtran}
\usepackage[cmex10]{amsmath}
\usepackage{multirow}
\usepackage{array}
\usepackage{sidecap}
\usepackage{url}

\newcommand{\bmm}[1]{ {\mbox{\boldmath $ {#1}$}}}

\newcommand{\bmath}[1]{ {\boldmath $ {#1}$}}

\newcommand{\beq}{\begin{equation}}
\newcommand{\eeq}{\end{equation}}
\usepackage{graphicx}        
\newcommand{\beqnr}{\begin{eqnarray}}
\newcommand{\eeqnr}{\end{eqnarray}}

\newcommand{\benum}{\begin{enumerate}}
\newcommand{\eenum}{\end{enumerate}}

\newcommand{\cW}{{\cal W}}

\newcommand{\qed}{\mbox{}\hspace*{\fill}\nolinebreak\mbox{$\rule{0.7em}{0.7em}$}}

\sidecaptionvpos{figure}{c}

\begin{document}
\title{Robust linear control of nonconvex battery operation in transmission systems}
\author{\IEEEauthorblockN{Daniel Bienstock, Gonzalo Mu\~noz, Shuoguang Yang}
\IEEEauthorblockA{Department of IEOR, Columbia University\\
New York, NY 10027\\
Email: $\{\mbox{dano,gm2543,sy2614}\}$@columbia.edu}
\and
\IEEEauthorblockN{Carsten Matke}
\IEEEauthorblockA{NEXT ENERGY\\
Oldenburg 26129 Germany\\
Email: carsten.matke@next-energy.de}
}
\maketitle

\begin{abstract} We describe a robust multiperiod transmission planning model
  including renewables and batteries, where battery output is used to
  partly offset renewable output deviations from forecast. A central element is
  a nonconvex battery operation model which is used with a robust
  model of forecast errors and a linear control scheme.  Even though the problem is nonconvex we provide an efficient and theoretically valid algorithm that
  effectively solves cases on large transmission systems.
\end{abstract}
\section{Introduction}
A great deal of recent work has focused on the integration of storage into
transmission systems, often in conjunction with the use of renewables.
A partial list includes \cite{Brown2008,Dufo2007,hedman,Jannati2014,guan,Siemer2016, Bludszuweit2011,Wu2014, dvorkin2, papavasiliou}.  The main goal of this work is 
to incorporate nonlinear and nonconvex battery operation models within a control framework that accounts for uncertainty in grid operation, specifically with regards to renewables. The control framework we study in this paper works as follows (with assumptions discussed below).
\begin{enumerate}
\item We consider a time horizon comprised by $T$ periods, each of length $\Delta$ (on the order of half an hour to one hour).  The output of each
generator as well as a linear control scheme governing each battery are
computed at time zero, using a forecast for renewable outputs.
\item At the start of each time period the average levels of renewable outputs
  for that period
  are estimated from real-time readings.  These estimations are used to set
  the output for each battery, through the control policy described above,
  and so as to offset the deviations of renewable outputs from the forecast. The
  output of each battery will be held constant during the current period, with
  additional real-time changes in renewable output handled through a
  standard scheme such as frequency control.
\end{enumerate}
This model simplifies many details, for brevity. For example, it assumes that loads are not subject to uncertainty.  Furthermore, we propose a
linear control, rather than affine.  An affine control governing battery operation would result on batteries supporting loads (and/or charged by standard
generators).  This extension is straightforward.  We omit these
extensions for brevity.

This broad line of control modeling is similar to recent work;   see e.g. \cite{andy, warrington, johanna}.  Our main contributions are as follows:
\begin{itemize}
\item [{\bf (a)}]  We  optimize over an explicitly nonconvex model for battery operation.  Batteries exhibit numerous complex nonlinear behaviors that
  can be very difficult to incorporate into an optimization framework. A standard model used in the literature relies on a constant charging and a constant
  discharging \textit{efficiency}.  If a (possibly negative) amount $E$ of energy is input into a
  battery, then the energy state of the battery changes by the amount
  $ \eta_c E^+ -\eta_d^{-1} E^{-}$ where $0 < \eta^c \le 1$ and
  (resp. $0 < \eta^d \le 1$) is the charging (discharging) efficiency and $E^+$
  and $E^-$ are the positive and negative parts of $E$.  Such a model is
  inherently nonconvex and its incorporation in optimization requires
  the complementarity condition $E^+ E^- = 0$.  This condition is sometimes modeled
  through the use of binary variables \cite{scuc} and is not a guaranteed outcome from
  a convex formulation.  We generalize this standard model by allowing state-dependent  charging and
  discharging efficiencies.  Additionally, we use a nonconvex charging \textit{speed} model.
\item [{\bf (b)}] We use a nonsymmetric and nonconvex robust model for renewable output
  deviations from forecast.  Here we note that popular stochastic distributions
  for wind power are nonsymmetric, e.g. Weibull distributions.  The use of
  nonsymmetries allows for specific risk stances with respect to renewable
  shortfall or excess.
\item [{\bf (c)}] Despite the above nonconvexities, we describe a theoretically
  valid and computationally practicable optimization scheme that reduces our
  overall problem into a sequence of convex, linearly
  constrained optimization problems.  Our algorithm is tested on realistic large
  transmission systems.
\end{itemize}
The detailed battery model is given in Section \ref{bmodel} and the
forecast errors model is described in Section \ref{dmodel}.
Section \ref{omodel} presents our optimization model, our algorithm
is given in Section \ref{algo} and experiments in Section \ref{exps}.

\section{Nomenclature}
\begin{tabular}{ll}
  $B$                 & bus susceptance matrix \\
  $P^g_{k,t} $       & standard generation  at bus $k$, period $t$\\
  $P^d_{k,t} $       & load  at bus $k$, period $t$\\
  $\theta^t_k$            & phase angle at bus $k$, period $t$\\
  $\bar w_{k,t} + $\bmm{w_{k,t}} & renewable output at bus $k$, period $t$:\\
  & $\bar w_{k,t}$ = forecast, \bmm{w_{k,t}} =  deviation\\
  $\cW$ &  set of deviations \bmath{w_{k,t}}\\
  $ L_{km}$          & line limit for line $km$\\
  $\lambda^t_{i,j}$     & battery control at buses $i, ~j$, period $t$
\end{tabular}  
\section{Battery model}\label{bmodel}
The accurate modeling of battery behavior and operating constraints is a nontrivial task,
made difficult by the wide variety of technologies and the complex nature
of the underlying chemical processes.  See e.g. \cite{winter}, \cite{Luo2015}.  Our
model seeks to incorporate relevant battery chemistry details
while resulting in a computationally accessible formulation.
We use charge-dependent piecewise-linear (or -constant) models for
charge/discharge \textit{efficiency} and
charge/discharge \textit{speed}. See e.g. \cite{powersonic} (page 19).\\

\begin{SCfigure}[][h!]
\centering
\includegraphics[scale=.33]{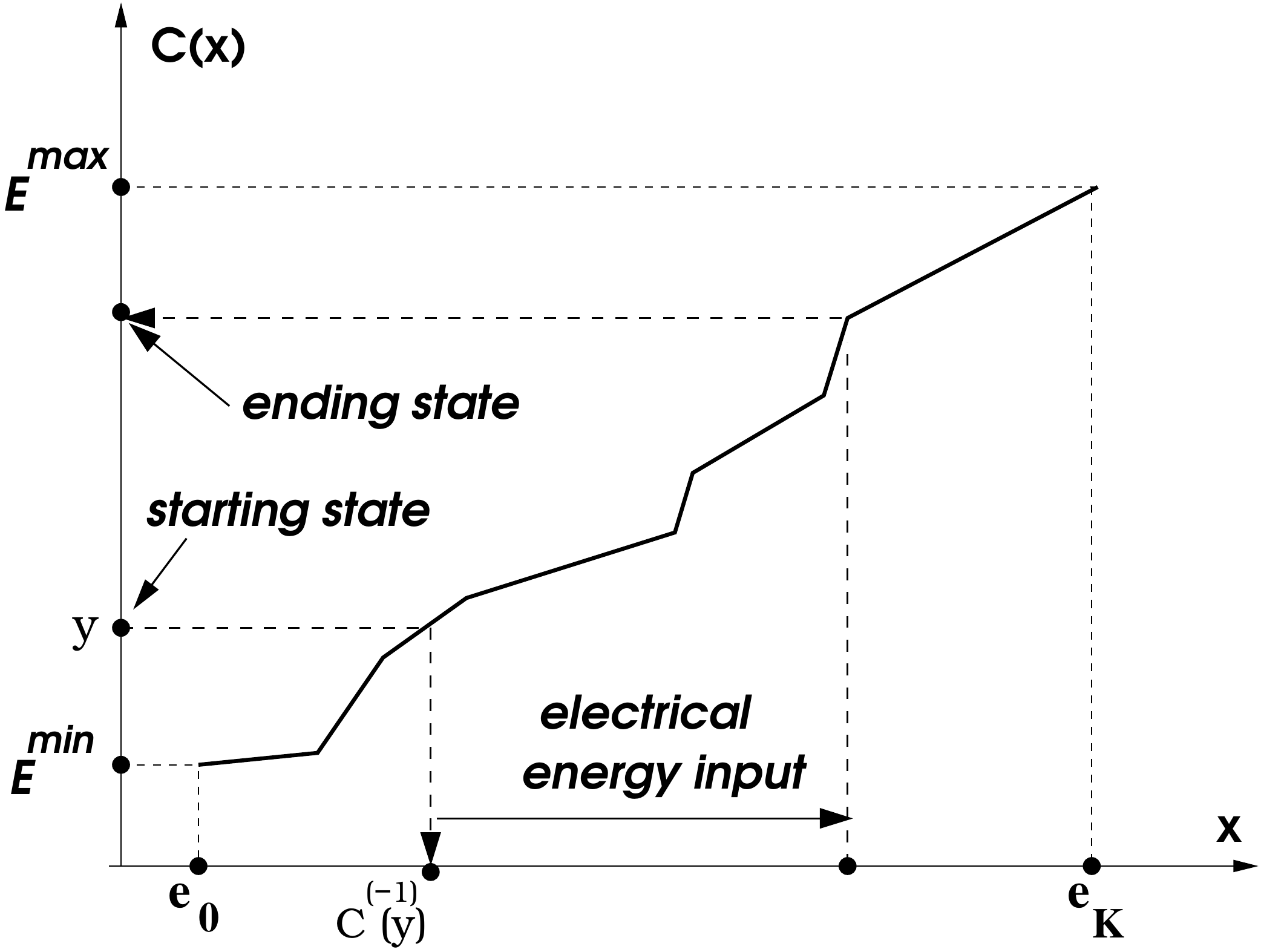}
\label{chargefig}
\end{SCfigure}

\noindent {\bf Charge efficiency.} Generalizing the standard model described above, we assume \textit{piecewise-constant} charging/discharging efficiencies.
Charging is described by a battery-specific monotonically increasing, piecewise linear battery {\bf charging function} $C(x)$, where $x$ represents electrical energy injection.    Suppose that the battery holds (chemical energy) charge $y$, and that we input electrical energy (= power $\times$ time) $\Sigma \ge 0$. Then
the charge of the battery will increase to $y + C( C^{(-1)}(y) + \Sigma)$, where
$C^{(-1)}$ is the inverse function of $C$.  See Figure \ref{chargefig}. 
The derivative of $C(x)$ (when it exists) is the charging efficiency at that point of the curve. We will assume a similar discharge function.

As in prior work (see e.g. \cite{scuc})
we are given an operating range $[ E^{\min}, E^{\max}]$ bor battery
charge.  Our detailed model works as follows. Let $[e_0, e_1, \ldots , e_K]$ where
$e_0 = C^{(-1)}(E^{\min})$ and $e_K = C^{(-1)}(E^{\max})$ be the breakpoints of
the charging function.  We impose conditions (a) and (b):\\

\noindent {\bf (a)} We require that the battery charge always remain
in the range $[E^{\min}, E^{\max}]$.\\

\noindent {\bf (b)} For each $t > 1$, suppose that at
the start of period $t$ battery charge lies in the range $[C(e_s), C(e_{s+1})]$ for some $s$.  Let $t'>t$,
and suppose our control mechanism  never discharges the battery in periods $t, t+1, \ldots, t'$. Then we require that the
maximum electrical energy input into the battery, in those periods, is at most
$C^{(-1)}(E^{\max}) - e_s$.  A similar statement is made in case we never
charge the battery in periods $t, t+1, \ldots, t'$.\\

Constraint (b) amounts to a more detailed approximation of the charge and
discharge functions than provided by (a) alone.\\

\noindent {\bf Charge speed.} Here we want to constrain the
maximum increase in battery charge in a period of length $\Delta$ as a function of
the battery charge at the start of the period.  It will be
computationally more convenient to state it in terms of input electrical
energy.   We impose:\\

\noindent {\bf (c)} Suppose that at the start of some period $t$  the charge of the battery, 
is in the range $[C(e_s) , C(e_{s+1})]$ (for some $s$). 
Then during period $t$ we can input into the battery electrical energy at most $v_s$. Here, 
$v_0, \ldots, v_K$ are given limits and for convenience $e_0, \ldots, v_K$
are the breakpoints for the charging function. The definition is necessarily ambiguous at the breakpoints, and it implies
state-dependent maximum (instantaneous) \textit{power injection} limits.
  
\section{Data robustness model} \label{dmodel}
Our linear control relies on \textit{estimations} on the quantities \bmath{w_{k,t}},
formally defined as follows.\\

\noindent {\bf Definition}. \bmath{w_{k,t}} is the average deviation of
renewable output at bus $k$ during period $t$.\\

\noindent We must therefore account for intrinsic stochastic variability on
renewable output and also on measurement errors and noise.  We rely
on a robust optimization model, which we term a \textit{concentration model},
which is given by {\bf nonnegative matrices} $K^+$ and $K^-$,
and a vector $b$, and
corresponds to the set $\cW$ of all $w$ satisfying
  \begin{eqnarray}
    &&    K^+ w^+ \ + \ K^- w^- \le \ b \label{notlinmodel}
  \end{eqnarray}
Here, $w^+$ is the vector with entries $w_{k,t}^+$
$= \max\{ w_{k,t} , 0 \}$ and likewise with $w^-$. We note that the description 
\eqref{notlinmodel} nonconvex. It can be used to model bounds on the individual
quantities $w^+_{k,t}$ and $w^-_{k,t}$, and allows for asymmetries
and correlation both across time and buses.   We assume that the set
of $w$ satisfying \eqref{notlinmodel} is \textit{full dimensional}. 

A special (though symmetric) case of this model is given by ``uncertainty budgets''
models (see \cite{Ben2009}, \cite{Bertsimas2011}), given for example
by the conditions
\begin{eqnarray}
   | w_{k,t}| & \le & \gamma_{t,k}, \quad  \mbox{all $t$ and $k$} \label{budg1}\\
  \sum_k (\gamma_{k,t})^{-1} |w_{k,t}| & \le & \Gamma^t \quad \mbox{all $t$} \label{budg2}
\end{eqnarray}
\noindent Here, the $\gamma^{k,t}$ and $\Gamma^t$ are parameters used to
model risk aversion, which can be chosen using data-driven techniques.
See \cite{datadriven}.  We can easily adapt these constraints so as to
include cross-time correlation.

Alternatively we could rely on \textit{chance-constrained} models (or, better,
on distributionally robust chance-constrained models).  A technical hazard
arising from the modeling of battery behavior can be outlined as follows.
Suppose that \bmath{x(t)}, $t = 1, \ldots, T$ is a discrete-time stochastic process.
Even if we
understand this process well enough so that e.g. we can compute tail probabilities for each individual \bmath{x(t)}, the \textit{partial sums }
\bmath{S(k) \ = \  \sum_{t = 1}^k x(t)} will in general be much more complex
random variables, except in special cases (i.e. gaussian distributions).  This
difficulty, combined with the goal of safe battery modeling, and the fact
that we need to account for observation errors, leads us to rely on a
robust model in this work.
\section{Optimization model} \label{omodel}
Our optimization model will compute  outputs $P^g_{k,t}$
for standard generation for each bus $k$ and period $t$ (fixed at zero if there
is no generator at that bus) and control parameters $\lambda^t_{ij}$ for
each period $t$ and pair of buses $i$, $j$.  This computation is assumed to take
place at time zero based on forecast data as discussed above. We next describe the generic
control scheme as well as a practical special case.

In the most general
case, our control works as follows: at the start of period $t$ we
compute
\textit{estimates} $w_{j,t}$
for all the values \bmath{w_{j,t}}, and during period $t$ the electrical power output of the
battery at bus $i$
will be set to \begin{eqnarray}
  && -\sum_{j} \lambda^t_{i,j} w_{j,t}, \label{lincontrol}
\end{eqnarray}  
  where $\lambda^t_{i,j} \ge 0$.  A simplified scheme relies, 
for each bus $i$ holding a battery, on a set $R(i)$, the set of renewable
buses that the battery at $i$ \textit{responds to}.  In this scheme 
the output at battery at $i$ in period $t$ will be of the form
$- \lambda^t_{i} \sum_{j \in R(i)} w_{j,t}$.  This simplified scheme only needs
an estimate of the (random) aggregate quantities $\sum_{j \in R(i)} \bmm{w_{j,t}}$.  Even more specialized cases are those where $R(i)$ =
all buses, and where $R(i) = i$ (a battery at each renewable bus). Regardless
of the special case, let us denote by $\Lambda^t$ the matrix with entries
$\lambda^t_{i,j}$, and by $w^t$ the vector with entries $w_{j,t}$.  In what follows, $\Lambda$ is a vector that includes all parameters $\lambda^t_{i,j}$. \\

Our
optimization problem is a robust multi-period DC-OPF-like problem, which
we term {\bf BATTOPF}:
\begin{subequations}\label{battopf}
\begin{eqnarray}
\hspace*{-.27in}  && \min_{P^g, \Lambda} \quad \sum_t \sum_k c_{k,t} (P^{g}_{k,t})  \label{battobj}\\
\hspace*{-.27in}    &&  \nonumber \\
\hspace*{-.27in}    &&  \mbox{s.t.} \quad 0 \, \le  P^g_{k,t} \, \le \, P^{g,\max}_{k,t} \quad \mbox{all period $t$ and buses $k$,} \label{battPbounds}\\
\hspace*{-.27in}    &&  0 \, \le \, \lambda^t_{i,j} \quad \mbox{all period $t$ and buses $i, \, j$,} \label{battlambdabounds}\\
\hspace*{-.27in}    &&  \mbox{and \eqref{stochbalance}-\eqref{stochbatt} feasible for all period $t$, and all $\bmm{w} \in \cW$:} \nonumber\\ 
\hspace*{-.27in}    && \hspace*{.2in} \bullet \ B \, \bmm{\theta}^t \ = \ P^{g}_t \, + \, \bar w_t + \bmm{w}_t - \Lambda^t \bmm{w}_{t} \ - \ P^{d}_t \quad \label{stochbalance} \\
\hspace*{-.27in}    && \hspace*{.2in} \bullet \ \mbox{line limits in period $t$, using phase angles $\theta^t$} \label{stochthermal} \\ 
\hspace*{-.27in}    && \hspace*{.2in} \bullet \ \mbox{battery operation constraints in period $t$} \label{stochbatt}
\end{eqnarray}
\end{subequations}
In this formulation the $c_{k,t}$ are standard OPF generation cost functions
(which we assume convex), $P^{g,\max}_{k,t}$ is the maximum generation at
bus $k$ in period $t$, $\Lambda^t$ is the matrix with entries $\lambda^t_{i,j}$, and $P^g_t, P^d_t$ and $\bmm{w_t}$ are (respectively) the vectors with entries $P^g_{t,k}, P^d_{t,k}$ and $\bmm{w_{k,t}}$.  The
formulation seeks a conventional generation plan plus a control algorithm
so as to minimize generation costs subject to remaining feasible  under
all data scenarios.  Constraint \eqref{stochbalance} is the standard DC flow
balance system.

Problem BATTOPF presents challenges because of the nonconvexities in the modeling of \eqref{stochbatt}.  Further, the concentration model $\cW$
is provided through a nonconvex description, and, as a result, the standard
linear programming duality approach in robust optimization cannot be applied.  In
Section \ref{algo} we will provide an efficient algorithm for solving this
problem.
We remark that \eqref{stochbalance}, which must hold for all $\bmm{w} \in \cW$, requires that for all $t$,
$$ 0 \ = \ \sum_{i} \left( P^{g}_{t,i} \, + \, \bar w_{t,i} + \bmm{w}_{t,i} - \sum_{j} \lambda^t_{i,j} \bmm{w}_{j,t} \ - \ P^{d}_{t,i} \right). $$
Since $0 \in \cW$ in particular we have that
$0 =  \sum_{i} \left( P^{g}_{t,i} \, + \, \bar w_{i,t} \ - \ P^{d}_{t,i} \right)$,
and as a result for all $\bmm{w} \in \cW$
\begin{eqnarray}
  0 & = & \sum_{i} \bmm{w}_{j,t}\left(1 - \sum_{i} \lambda^t_{i,j} \right). \label{reduced}
\end{eqnarray}
Since the
set $\cW$ is assumed full-dimensional,
and \eqref{reduced} describes a hyperplane in $\bmm{w}$-space,
we conclude that \eqref{stochbalance}
is feasible (for a given $P^g$ and $\Lambda$) if and only if
\begin{eqnarray}
  1 & = &  \sum_{i} \lambda^t_{i,j} \quad \forall \ t, \ j.\label{red2}
\end{eqnarray}

\section{Algorithm}\label{algo}
As described above it appears difficult to produce an explicit, practicable
convex formulation for BATTOPF.  Here instead we provide
an efficient cutting-plane procedure.  Even though the outline of the
algorithm below is standard \cite{benders}, Steps {\bf 2} and {\bf 3} are
novel and critical.

The algorithm relies on a linearly constrained relaxation for BATTOPF
termed the \textit{master formulation}, with
objective \eqref{battobj}. At the start of the
procedure the master formulation includes constraints \eqref{battPbounds},
\eqref{battlambdabounds} and \eqref{red2}.  Each iteration produces a
\textit{candidate solution} $(\tilde P^g, \tilde \Lambda)$ for BATTOPF.
If this candidate is \textit{infeasible} for BATTOPF,
that is to say there is a realization of renewable deviations under which
the candidate
fails to satisfy \eqref{stochthermal}-\eqref{stochbatt}, then the algorithm
identifies a
linear inequality that is valid for all feasible solutions for BATTOPF,
yet violated by the candidate.  This inequality is then added
to the master formulation.  Thus, at each iteration the master formulation is
a relaxation for BATTOPF, and hence if at some iteration the
current candidate is feasible for BATTOPF
then it is optimal.  Formally, the algorithm is as follows:\\

\noindent {\bf 1.} Solve the master formulation, with solution $(\tilde P^g, \tilde \Lambda)$.

\noindent {\bf 2.} Check whether $(\tilde P^g, \tilde \Lambda)$ is feasible
for BATTOPF, that is to say, $(\tilde P^g, \tilde \Lambda)$
satisfies \eqref{stochthermal}-\eqref{stochbatt} for every $w \in \cW$.  If
so $(\tilde P^g, \tilde \Lambda)$ is {\bf optimal} for BATTOPF.
{\bf STOP}.

\noindent {\bf 3.}  Otherwise, there is $\hat w \in \cW$ such that
$(\tilde P^g, \tilde \Lambda)$ does not satisfy \eqref{stochthermal}-\eqref{stochbatt} when the renewable deviations are given by $\hat w$.  Compute an inequality
\begin{eqnarray}
&& \alpha^T P^g \ + \ \beta^T \Lambda \ \ge \ \alpha_0 \label{cut}
\end{eqnarray}
which is satisfied by all solutions to BATTOPF, but
violated by
$(\tilde P^g, \tilde \Lambda)$, add it to the master formulation, and go to {\bf 1}. \\

Steps {\bf 2} and {\bf 3} are both nontrivial because of the nonconvexity of the
battery model and the nonconvex description of the uncertainty set $\cW$. See
Section \ref{tech} for details.

\section{Numerical Experiments} \label{exps}
The algorithm was implemented using Gurobi \cite{gurobi} as
the LP solver. The first set of numerical tests study the scalability of the
algorithm as the number of time periods increases.  For these tests we used
the winter peak Polish grid from MATPOWER \cite{MATPOWER},
with $2746$ buses, $3514$ branches, $388$ generators,
 base load
(approx.) $24.8$ GW, 32 wind farms, with forecast output $ 4.5$ GW and 32
batteries, with total initial charge approx. $3.2$ GWh after unit conversion.  We used the uncertainty
budgets robustness model with an implied forecast error of up to $8.9 \%$. In
the runs below, loads increase (in similar proportions) in the first six periods.  For these runs we used one-piece charging/discharging curves. In our implementation, the initial formulation includes all line constraints
and equations \eqref{stochbalance}
for the nominal case (no forecast errors). In the following table, ``n'' and ``m'' indicated the number of variables and
constraints in the master formulation at termination (but including some preprocessing).
\begin{table}[htb]
%
  %
  \centering
\begin{tabular}{c|c|c|c|c|c}
\hline\noalign{\smallskip}
 {\bf T}    &{\bf n} & {\bf m} &{\bf Cost} &  {\bf Iterations} &  {\bf Time (s)}\\ 
\noalign{\smallskip}\hline\noalign{\smallskip}
6 & 20940 & 57519 &  7263172 &  20 &  366 \\
8 & 27920  & 76651 & 9738804  &  17 & 442  \\
10 & 34900 & 95825 &  12260289  & 19  & 670  \\
12 & 41880 & 114995 & 14784028  & 18  & 848  \\
\noalign{\smallskip} \hline\noalign{\smallskip}
\end{tabular}
\end{table}

The running time is primarily accrued by solving linear programs, whose size
grows proportional to the number of periods.  The number of iterations
appears nearly constant.

Next we describe, in greater detail, a one-period example derived from the ``Case9''
dataset in Matpower.  This case was modified by adding renewables and
batteries.
\begin{itemize}
  \item Renewables are located at buses 4 and 8,  with forecaset output
  $50$ and $100$ MW, respectively.  We using a concentration model
  given by
  the constraints $-50 \le \bmm{w_{4,1}} \le 0$,
  $-100 \le \bmm{w_{5,1}} \le 0$, and $ 2 |\bmm{w_{4,1}}| + |\bmm{w_{5,1}}| \ \le \ 100$. In particular the system may experience a loss of up to $100$ MW in renewable power.
\item Identical batteries are located at buses 4 and 9, with charging efficiency $1.0$ and discharging efficiency $0.8$ (note that since we only consider
  renewable decreases, only discharges will take place). Both batteries start
  with $80$ units of charge, and the instantaneous maximum (power) discharge rate
  is set at $100$ MW.  Each of the batteries to the aggregation of renewable
  errors, i.e. the output of each battery $i$ is of the form $-\lambda_{i} (\bmm{w_{4,1}} + \bmm{w_{5,1}})$.
\item The limits of lines $4-5$,  $5-6$, $6-7$,  $7-8$, $8-9$ and $9-4$ were
  reduced to $50, 75, 50, 90, 100$ and $70$ MW (resp.)
\end{itemize}
The forecast renewable generation, $150$ MW, amounts to almost $50 \%$ of
total load ($315$ MW).  The minimum (DC-OPF) generation cost for Case9 is
approximately $5216$, and due to the large renewable penetration, in our example
the cost in the nominal case (no forecast errors) is much lower:  $2384.75$.
However, this solution is certainly not robust.

The solution to the robust problem has cost $2488.05$ (less than one percent
increase over the non-robust solution) and is given by
$\lambda_{4} = 0.36$, $\lambda_{9} = 0.64$, $P^g$.  We can briefly examine the
feasibility of this solution as follows. Consider, first, the battery charge
constraints, and suppose renewable output drops by
$100$ MW (the maximum allowed by the model), then the charge of the battery at bus 9 will decrease by
$0.64 \times 100/.8 = 80$ units. Since the battery starts with $80$ units it will
therefore drain completely (but not go negative). Likewise, the charge
of the battery at bus $4$ will likewise drop by $0.36 \times 100/.8 = 45 $ units from its original $80$ units.

\section{Technical details}\label{tech} Here 
we outline efficient procedures for {\bf 2} and {\bf 3}.  The battery
constraints are the most delicate. For brevity we focus on battery speed, and the next result is key. \\

\noindent {\bf Lemma 2.}  Suppose that $(\tilde P^g, \tilde \Lambda)$ is
an infeasible candidate solution for BATTOPF.  Let $\hat w
\in \cW$ be a realization of renewable deviations under which
$(\tilde P^g, \tilde \Lambda)$ fails to satisfy the battery operation
constraint at some period $t$ and bus $k$.  Then, without loss of generality, $\hat w$ is {\bf sign consistent} for all periods $h < t$,
that is to say, either $\hat w_{i,h} \ge 0$ for all buses $i$ and periods $h < t$,
or $\hat w_{i,h} \le 0$ for all buses $i$ and periods $h < t$.  \\

\noindent Lemma 2 (proved later) has an important corollary. Namely, Step {\bf 2} gives rise,
for each time period, to
two linear programs for each battery speed constraint. We discuss these points next.

Consider, first, the battery speed constraint for a battery at bus $i$ and period $t$.  Given $\tilde \Lambda$ this constraint will fail to hold if there is
$\hat w \in \cW$ such that in period $t$ the battery charge is in some interval
$(C(e_s), C(e_{s+1}))$ and yet the electrical power input into the battery,
as per $\tilde \Lambda$, exceeds
$v_s$.  Since
we can assume that $\hat w$ is sign-consistent consider the case where all
$w_{j,i} \ge 0$ for $i < t$.  If at time zero the battery has charge $E_{0}$ then
\begin{subequations} \label{thefoo}
\begin{eqnarray}
\hspace*{-.4in}   && e_s \ \le \ E_{0} \ + \ \sum_{h = 1}^{t-1} \sum_{j} \tilde \lambda^t_{i,j} \hat w_{j,t} \ \le \  e_{s+1} \quad \mbox{and} \label{foobar}\\
\hspace*{-.4in}   && \sum_{j} \lambda^t_{i,j} \hat w_{j,t} \ > \ v_s. 
\end{eqnarray}
\end{subequations}
Thus, the existence of such a vector $\hat w$ can be tested by solving the linear
program where we maximize $\sum_{j} \lambda^t_{i,j} w_{j,t}$ subject to \eqref{foobar}, $ K^+ \hat w \, \le \, b$ (the positive part of the concentration model), and $\hat w \ge 0$.  The case $\hat w \le 0$ is
similar.  If a vector $\hat w$ satisfying \eqref{thefoo} is found, then
it proves that any control $\Lambda$ feasible for BATTOPF must violate at least one of the three inequalities in \eqref{thefoo}, giving rise to a
``disjunctive cut'' \cite{balas}.\\

In terms of Steps 2 and 3
of our algorithm
we have to consider, lastly, line limits.  To fix ideas consider the constraint that the power flow on line $km$ does
not exceed its limit $L_{km}$ in period $t$.  Using shift factors (i.e. a pseudoinverse for $B$) we need
to check that
\begin{eqnarray}
  \nu_{km}^T \left(  \tilde P^{g}_t \, - \,  P^d_t  + \, {\bar w_t} - \tilde \Lambda_t w_t \right) & \le & L_{km} \label{pkm}
\end{eqnarray}  
for all $w \in \cW$, where $\nu_{km}$ is an appropriate vector.  It
sufficies to maximize the right-hand side of \eqref{pkm} over all $w \in \cW$.  This can be done by solving the LP
\begin{eqnarray}
\hspace*{-.2in}  && \max \quad \nu_{km} \left(  \tilde P^{g}_t \, - \,  P^d_t  + \, {\bar w_t} - \tilde \Lambda_t \hat w_t \right) \nonumber \\
\hspace*{-.2in}    &&  \mbox{s.t.} \quad \quad K^+ \hat w^{(p)} \ + \ K^- \hat w^{(n)} \le \ b \label{notlinmodel2}\\
\hspace*{-.2in}    &&  \quad \quad \quad \hat w  \ = \ \hat w^{(p)} - \hat w^{(n)}, \label{bal} \\
\hspace*{-.2in}    &&  \quad \quad \quad \hat w^{(p)} \ge 0, \quad  \hat w^{(n)} \ge 0.
\end{eqnarray}
Here $\hat w$ represents $w$, and $\hat w^{(p)}$ and
variables 
$\hat w^{(n)}$ represent, respectively, $w^+$ and $w^-$.  Since 
$K^+ \ge 0$ and $K^- \ge 0$ there is an optimal solution to this
LP satisfying the complementarity condition
$\hat w^{(p)}_{k,h} \hat w^{(n)}_{k,h} = 0$ for all buses $k$ and periods $h$.
Thus $\hat w \in \cW$ and \eqref{notlinmodel2} yields the concentration model
\eqref{notlinmodel}.

If the value of the LP exceeds $L_{km}$, then, where $\hat w \in \cW$ is
an optimal solution satisfying complementarity, 
the inequality 
\begin{eqnarray}
\nu_{km} \left(  \tilde P^{g}_t \, - \,  P^d_t  + \, {\bar w_t} - \tilde \Lambda_t \tilde w_t \right) & \le &   L_{km} \label{pkm2}
\end{eqnarray}  
is valid for BATTOPF but violated
by $(\tilde P^{g}, \tilde \Lambda)$.\\

\noindent {\bf Proof of Lemma 2.}  We sketch a proof of Lemma 2 in the case of
the charge speed constraints.  Suppose that $\tilde \Lambda$ violates the charge speed constraint for some battery $k$ during period $t$, under deviations $w \in \cW$.
Let $E_0$ be the battery charge at time zero (known), and the energy state
at the start of period $t$, under deviations $w$, be $E'$. Then
$E' \in (e_s, e_{s+1})$ (for some $s$) and the energy input into the battery
during period $t$ exceeds the maximum $v_s$. Assume  that $C^{(-1)}(E') \ge C^{(-1)}(E_0)$ (the reverse case is similar).  Let $\bar w_{i,h} = w_{i,h}^+$ for all buses $i$ and periods $h < t$, and $\bar w_{i,h} = w_{i,h}$
  otherwise.  Then $\bar w \in \cW$.  
  The deviations $\bar w$ do not cause discharging before period $t$, and 
  so the at the start of period $t$ the energy state will be at least
$E'$.  Since the charging function is monotonically increasing and continuous, for some value $0 \le \kappa \le 1$ the vector $\breve w$ defined
by $\breve w_{i,h}  = \kappa \, \bar w_{i,h}$ for $h < t$, and
$\breve w_{i,h} = w_{i,h}$ otherwise, is such that 
that under deviations $\breve w$ the energy state at the start of period
$t$ will be exactly $E'$. Yet, clearly, $\breve w \in \cW$, $\breve w_{i,h} \ge 0$ for all $i$ and $h < t$.

\bibliographystyle{IEEEtran}
\bibliography{IEEEabrv,batt}

\end{document}